\documentclass[12pt]{amsart}
\makeatletter
\let\atopwithdelims\@@atopwithdelims
\let\over\@@over

\newtheorem{theorem}{Theorem}
\newtheorem{lemma}[theorem]{Lemma}
\newtheorem{corollary}[theorem]{Corollary}

\newenvironment{example}{\mbox{\sc Example} }{}

\newcommand{\fv}{{{\bf f}}}
\newcommand{\gv}{{{\bf g}}}

\newcommand{\Z }{{Z\!\!\! Z}}
\newcommand{\pa }{{\partial}}

\newcommand{\lb }{{\lambda}}
\newcommand{\da }{{\delta}}

\begin{document}

\title[The number of $k$-faces of a simple $d$-polytope]{The number of
{\large $\lowercase{k}$}-faces of a simple {\large $\lowercase{d}$}-polytope}
\author{Anders Bj\"orner}
\address{\hskip-\parindent Anders Bj\"orner\\
Department of Mathematics \\
Royal Institute of Technology \\
S-100 44 Stockholm, SWEDEN}
\email {bjorner@math.kth.se }
\author{Svante Linusson}
\address{\hskip-\parindent Svante Linusson\\
Department of Mathematics \\ 
Stockholm University \\ 
S-106 91
Stockholm, SWEDEN}
\email{linusson@matematik.su.se}

\thanks{Both authors were supported by EC grant CHRX-CT93-0400
and by the Mathematical Sciences Research Institute (Berkeley, CA).
Research at MSRI is supported in part by NSF grant DMS-9022140.}

\begin{abstract}
Consider the question: Given integers $k<d<n$, does there exist a
simple $d$-polytope with $n$ faces of dimension $k$? We show that
there exist numbers $G(d,k)$ and $N(d,k)$ such that for $n> N(d,k)$
the answer is yes if and only if $n\equiv 0\quad  \pmod {G(d,k)}$.
Furthermore, a formula for $G(d,k)$ is given, showing that e.g. 
$G(d,k)=1$ if 
$k\ge \left\lfloor\frac{d+1}{2}\right\rfloor$ or if both $d$ and $k$
are even, and also in some other cases (meaning that all numbers beyond
$N(d,k)$ occur as the number of $k$-faces of some simple $d$-polytope).

This question has previously been studied only for the case 
of vertices ($k=0$), where Lee \cite{Le} proved the existence
of $N(d,0)$ (with $G(d,0)=1$ or $2$ depending on whether $d$ is
even or odd), and Prabhu
\cite{P2} showed that 
$N(d,0) \le cd\sqrt {d}$. We show here that asymptotically
the true value of Prabhu's  constant is $c=\sqrt2$ if $d$ is
even, and $c=1$ if $d$ is odd.
\end{abstract}

\maketitle

\section{Introduction}

An integer $n$  will be called $(d,k)$-{\it realizable} if there is a simple
$d$-polytope with $n$ faces of dimension $k$. 
For terminology and basic properties of polytopes we refer to the
literature, see e.g. \cite{Z}.

We show, see Theorem \ref{gap-thm}, that there 
exist numbers $G(d,k)$ and $N(d,k)$ such that
\begin{itemize}
\item{if $n$ is $(d,k)$-realizable then $G(d,k)$ divides $n$;}
\item{if $G(d,k)$ divides $n$ and $n>N(d,k)$ then $n$ is $(d,k)$-realizable.}
\end{itemize}
The $G(d,k)$-divisible numbers that are not $(d,k)$-realizable will 
be called $(d,k)$-{\it gaps}. Thus there are only finitely many gaps
for all $d>k\ge 0$. In this paper we study the numbers  $G(d,k)$ and $N(d,k)$.
Our proofs rely on the $g$-theorem.

To give some feeling for the results, let us discuss a few special cases.
The parity restrictions that exist for each dimension $k$ are easiest to
understand for the case of vertices ($k=0$). Namely, the graph of a simple
$d$-polytope is $d$-regular, so if the polytope has $n$ vertices then it
has $\frac{dn}{2}$ edges. Hence, if $d$ is odd $n$ must be even. This
is in fact the only constraint, and we have
\[   G(d,0)=\left\{
\begin{array}{rl}
 1, &\quad \mbox{$d$ even}\\
 2, &\quad\mbox{$d$ odd.}
\end{array}\right.
\]
This result is due to Lee \cite{Le}, who initiated the study of
properties of vertex-count numbers of simple polytopes. Via the regular graph
property  this also implies the result for edge-count numbers:
\[   G(d,1)=\left\{
\begin{array}{rl}
 \frac{d}{2}, &\quad \mbox{$d$ even}\\
 d, &\quad\mbox{$d$ odd.}
\end{array}\right.
\]

For $1<k<\left\lfloor\frac{d+1}{2}\right\rfloor$ the
situation gets more complicated and the answer is different for $k$ even
and $k$ odd. For instance,
\[   G(d,2)=\left\{
\begin{array}{rl}
 2, &\quad \mbox{$d\equiv 1 \pmod {4}$}\\
 1, &\quad\mbox{otherwise.}
\end{array}\right.
\]
The modulus $G(d,k)$ can get arbitrarily large in this range;
for instance, $G(d,k)=d-k+1$ whenever $k$ is odd and $d-k+1$ is a prime.
Then, for $k\ge \left\lfloor\frac{d+1}{2}\right\rfloor$ the situation 
simplifies again to $G(d,k)=1$.  Theorem \ref{th:main} gives the 
general formula for $G(d,k)$.

\smallskip
It is also of interest to study the magnitude of the numbers $N(d,k)$
(defined as the smallest possible ones for which the above statement is
true). Again, this has been studied for the case of vertices by
Prabhu \cite{P2}, who showed that 
$N(d,0) \le cd\sqrt {d}$. We prove that asymptotically the true value of
Prabhu's  constant is $c=\sqrt2$ if $d$ is even, and $c=1$ if $d$ is odd,
see Section \ref {sc:vertices}.
We also give an upper bound for $N(d,k)$ in the general case, 
Theorems \ref{N-thm} and \ref{th:stora_k}, but leave
open the determination of its true asymptotic growth.

\section{Preliminaries}
Given a $d$-dimensional polytope $P$, we will call 
${\fv}:=(f_0,f_1,\dots,f_{d-1})$ the $f$-{\em vector} of $P$,
where  $f_i$ is the number of faces of dimension $i$.

For any integers $n,s\ge 1$, there is a unique way of writing
\[ n={a_s\choose s}+{a_{s-1}\choose s-1}+\dots+{a_i\choose i},
\]
so that $a_s>a_{s-1}>\dots>a_i\ge i\ge 1$.  Then define:
\[\pa^s(n):={a_s-1\choose s-1}+{a_{s-1}-1\choose s-2}+
\dots+{a_i-1\choose i}.
\]
Also let $\pa^s(0):=0$.

A nonnegative integer sequence
$(n_0,n_1,n_2,\dots)$ is called an {\it $M$-sequence} if
\[n_0=1 \quad\mbox{ and } \quad \partial^s(n_s)\le n_{s-1}\quad \mbox{for
all }s\ge 1.
\]
Two simple facts we will need about $M$-sequences is that if there is a
zero in the sequence then all the following entries are also zeros, and
that any sequence satisfying $n_0=1$ and $n_1\ge n_2\ge n_3 \ge\dots$
is an $M$-sequence.

An alternative definition of $M$-sequence, due to Macaulay and
Stanley \cite{S1}, says that a sequence is a $M$-sequence if and only
if it is the $f$-vector of a multicomplex.  See \cite {Li} and \cite {Z}
for examples of other interpretations of $M$-sequences.
Let $\left\lfloor x\right\rfloor$ and $\left\lceil x\right\rceil$ 
denote the largest integer less than or equal to $x$ and the smallest 
integer larger than or equal to $x$, respectively.

Let $\da:=\left\lfloor {d \over 2}\right\rfloor$ and
let $M_d=(m_{ik})$ be the $(\da +1)\times d$-matrix with entries
\[ m_{ik}={d+1-i\choose k+1}-{i\choose k+1}, \quad \mbox{for }
0\le i\le \da,\ 0\le k \le d-1.
\] 
For example,
\[
M_{10}=\left( \begin{array}{rrrrrrrrrr}
11&55 &165 &330 &462 &462 & 330&165 &55 &11 \\
 9& 45& 120 &210 &252 &210 &120 &45 &10& 1 \\
 7&  35 & 84 &126 &126 & 84 &36 & 9 & 1 & 0 \\
 5  &25 & 55 & 70 & 56 &28 & 8 & 1 &0&0 \\
 3 & 15 & 31 & 34 &21 & 7 & 1 &0 &  0 &0\\
 1 &  5 & 10  &10 & 5 & 1 &0&0 &  0 &  0  \\
\end{array}
\right)
\]

Our proofs will rely on the $g$-theorem, conjectured by McMullen.
Sufficiency was proved by Billera and Lee \cite{BL}, and necessity by 
Stanley \cite{S2} and later by McMullen \cite{M}, 
see \cite{Z}.  We will use the
following matrix reformulation of the $g$-theorem, given by Bj\"orner
\cite{B1,B2}, see also \cite{Z}.  We have here reformulated the 
statement from simplicial polytopes to simple polytopes, which
just correspond to reading the $f$-vector backwards.

\medskip
\noindent

\begin{theorem}\label{g-thm}{\bf [The $g$-theorem]}
The matrix equation
\[ \fv=\gv\cdot M_d
\]
gives a one-to-one correspondence between $f$-vectors $\fv$ of simple
$d$-polytopes and $M$-sequences $\gv=(g_0,g_1,\dots,g_\da)$.\qed
\end{theorem}

\section{The modulus $G(d,k)$} 

The modulus mentioned in the introduction is defined as follows
\begin{equation}\label{Gdef}
G(d,k):=\gcd(m_{1,k},m_{2,k},\dots,m_{\da, k}),
\end{equation}
the greatest common divisor for
the elements in the $k$-th column and below the top row of the matrix $M_d$.
In this section we will give simple and explicit formulas for 
$G(d,k)$. The role of $G(d,k)$ 
as the period for the possible numbers of $k$-faces 
of $d$-polytopes will be shown in the next section.

\begin{theorem}\label{th:main} 
\begin{enumerate}
\item[(i)]{If $k\ge\left\lfloor\frac{d+1}{2}\right\rfloor$, then
$G(d,k)=1.$}
\item[(ii)]{If $k<\left\lfloor\frac{d+1}{2}\right\rfloor$ is
even, let $e$ be the  integer such that $2^{e}\le k+1<2^{e+1}$. Then 
\[G(d,k)=\left \{ 
\begin{array}{rl} 2, & \quad \mbox{ if $d-k+1\equiv 0$ (mod $2^{e+1}$)}\\ 1, &
\quad \mbox{ otherwise}.
\end{array}
\right.
\] }
\item[(iii)]{If $k<\left\lfloor\frac{d+1}{2}\right\rfloor$ is odd,  let
$p_1,\dots,p_t$ be the primes smaller than or equal to $k+1$, 
and let $e_i\ge 1$ be the integers such that 
$p_i^{e_i}\le k+1<p_i^{e_i+1}$. Then 
\[G(d,k)=
\frac{d-k+1}{\gcd(d-k+1, p_1^{e_1}p_2^{e_2}\cdots p_t^{e_t})}.
\] }
\end{enumerate}
\end{theorem}

For the proof we  will need some facts about binomial coefficients 
modulo powers of a prime, that will be developed in a sequence of lemmas. 

\begin{lemma}\label{lm:period} 
Let $k,e\ge 0$ and let $p$ be a prime 
such that $p^e\le k+1<p^{e+1}$. Then for all $r\ge 1$ we have
\begin{enumerate}
\item[(i)]{$d\equiv d' \;\pmod {p^{e+r}} \quad$ implies 
$\quad {d\choose k+1}\equiv {d'\choose k+1}\; \pmod {p^r};$
}
\item[(ii)]{$${p^{e+r}+k-i\choose k+1}\equiv (-1)^{k+1}{i\choose k+1}\; 
\pmod {p^r},$$
\newline
for all $i=0,1,\dots ,p^{e+r}+k$; }
\item[(iii)] the unique longest run of zeros in the period
is ${d\choose k+1}\equiv_{p^r} 0$ 
\newline for all 
$p^{e+r}\le d \le p^{e+r}+k$.\end{enumerate}
\end{lemma}

Part $(ii)$ shows that if $k$ is odd then the period extended by $k$
is symmetric $\pmod {p^{r}}$,
and if $k$ is even then it is antisymmetric. The lemma is
illustrated by the following modular Pascal triangle.

\begin{center}
\setlength{\unitlength}{1.3pt}
\begin{picture}(310,230)(-10,-10)
\put(135,223){1}
\put(130,215){1\ \ 1}
\put(125,207){1\ \ 2\ \ 1}
\put(120,199){1\ \ 3\ \ 3\ \ 1}
\put(115,191){1\ \ 0\ \ 2\ \ 0\ \ 1}
\put(110,183){1\ \ 1\ \ 2\ \ 2\ \ 1\ \ 1}
\put(105,175){1\ \ 2\ \ 3\ \ 0\ \ 3\ \ 2\ \ 1}
\put(100,167){1\ \ 3\ \ 1\ \ 3\ \ 3\ \ 1\ \ 3\ \ 1}
\put(95,159){1\ \ 0\ \ 0\ \ 0\ \ 2\ \ 0\ \ 0\ \ 0\ \ 1}
\put(90,151){1\ \ 1\ \ 0\ \ 0\ \ 2\ \ 2\ \ 0\ \ 0\ \ 1\ \ 1}
\put(85,143){1\ \ 2\ \ 1\ \ 0\ \ 2\ \ 0\ \ 2\ \ 0\ \ 1\ \ 2\ \ 1}
\put(80,135){1\ \ 3\ \ 3\ \ 1\ \ 2\ \ 2\ \ 2\ \ 2\ \ 1\ \ 3\ \ 3\ \ 1}
\put(75,127){1\ \ 0\ \ 2\ \ 0\ \ 3\ \ 0\ \ 0\ \ 0\ \ 3\ \ 0\ \ 2\ \ 0\
\ 1}
\put(70,119){1\ \ 1\ \ 2\ \ 2\ \ 3\ \ 3\ \ 0\ \ 0\ \ 3\ \ 3\ \ 2\ \ 2\
\ 1\ \ 1}
\put(65,111){1\ \ 2\ \ 3\ \ 0\ \ 1\ \ 2\ \ 3\ \ 0\ \ 3\ \ 2\ \ 1\ \ 0\
\ 3\ \ 2\ \ 1}
\put(60,103){1\ \ 3\ \ 1\ \ 3\ \ 1\ \ 3\ \ 1\ \ 3\ \ 3\ \ 1\ \ 3\ \ 1\
\ 3\ \ 1\ \ 3\ \ 1}
\put(55,95){1\ \ 0\ \ 0\ \ 0\ \ 0\ \ 0\ \ 0\ \ 0\ \ 2\ \ 0\ \ 0\ \ 0\
\ 0\ \ 0\ \ 0\ \ 0\ \ 1}
\put(50,87){1\ \ 1\ \ 0\ \ 0\ \ 0\ \ 0\ \ 0\ \ 0\ \ 2\ \ 2\ \ 0\ \ 0\
\ 0\ \ 0\ \ 0\ \ 0\ \ 1\ \ 1}
\put(45,79){1\ \ 2\ \ 1\ \ 0\ \ 0\ \ 0\ \ 0\ \ 0\ \ 2\ \ 0\ \ 2\ \ 0\
\ 0\ \ 0\ \ 0\ \ 0\ \ 1\ \ 2\ \ 1}
\put(40,71){1\ \ 3\ \ 3\ \ 1\ \ 0\ \ 0\ \ 0\ \ 0\ \ 2\ \ 2\ \ 2\ \ 2\
\ 0\ \ 0\ \ 0\ \ 0\ \ 1\ \ 3\ \ 3\ \ 1}
\put(35,63){1\ \ 0\ \ 2\ \ 0\ \ 1\ \ 0\ \ 0\ \ 0\ \ 2\ \ 0\ \ 0\ \ 0\
\ 2\ \ 0\ \ 0\ \ 0\ \ 1\ \ 0\ \ 2\ \ 0\ \ 1}
\put(30,55){1\ \ 1\ \ 2\ \ 2\ \ 1\ \ 1\ \ 0\ \ 0\ \ 2\ \ 2\ \ 0\ \ 0\
\ 2\ \ 2\ \ 0\ \ 0\ \ 1\ \ 1\ \ 2\ \ 2\ \ 1\ \ 1}
\put(25,47){1\ \ 2\ \ 3\ \ 0\ \ 3\ \ 2\ \ 1\ \ 0\ \ 2\ \ 0\ \ 2\ \ 0\
\ 2\ \ 0\ \ 2\ \ 0\ \ 1\ \ 2\ \ 3\ \ 0\ \ 3\ \ 2\ \ 1}
\put(20,39){1\ \ 3\ \ 1\ \ 3\ \ 3\ \ 1\ \ 3\ \ 1\ \ 2\ \ 2\ \ 2\ \ 2\
\ 2\ \ 2\ \ 2\ \ 2\ \ 1\ \ 3\ \ 1\ \ 3\ \ 3\ \ 1\ \ 3\ \ 1}
\put(15,31){1\ \ 0\ \ 0\ \ 0\ \ 2\ \ 0\ \ 0\ \ 0\ \ 3\ \ 0\ \ 0\ \ 0\
\ 0\ \ 0\ \ 0\ \ 0\ \ 3\ \ 0\ \ 0\ \ 0\ \ 2\ \ 0\ \ 0\ \ 0\ \ 1}
\put(10,23){1\ \ 1\ \ 0\ \ 0\ \ 2\ \ 2\ \ 0\ \ 0\ \ 3\ \ 3\ \ 0\ \ 0\
\ 0\ \ 0\ \ 0\ \ 0\ \ 3\ \ 3\ \ 0\ \ 0\ \ 2\ \ 2\ \ 0\ \ 0\ \ 1\ \ 1}
\put(5,15){1\ \ 2\ \ 1\ \ 0\ \ 2\ \ 0\ \ 2\ \ 0\ \ 3\ \ 2\ \ 3\ \ 0\
\ 0\ \ 0\ \ 0\ \ 0\ \ 3\ \ 2\ \ 3\ \ 0\ \ 2\ \ 0\ \ 2\ \ 0\ \ 1\ \ 2\
\ 1}
\put(0,7){1\ \ 3\ \ 3\ \ 1\ \ 2\ \ 2\ \ 2\ \ 2\ \ 3\ \ 1\ \ 1\ \ 3\
\ 0\ \ 0\ \ 0\ \ 0\ \ 3\ \ 1\ \ 1\ \ 3\ \ 2\ \ 2\ \ 2\ \ 2\ \ 1\ \ 3\
\ 3\ \ 1}
\put(-5,-1){1\ \ 0\ \ 2\ \ 0\ \ 3\ \ 0\ \ 0\ \ 0\ \ 1\ \ 0\ \ 2\ \ 0
\ 3\ \ 0\ \ 0\ \ 0\ \ 3\ \ 0\ \ 2\ \ 0\ \ 1\ \ 0\ \ 0\ \ 0\ \ 3\ \ 0\
\ 2\ \ 0\ \ 1}
\end{picture} 
\end{center}
\begin{center}
Pascal's triangle (mod 4).
\end{center}

\bigskip
For each prime $p$  define the valuation 
$v_p:\Z\setminus\{0\}\rightarrow\mathbb N$ by $v_p(n)=s$, where $p^s$ is the
highest power of $p$ that is a divisor to
$n$.  We will frequently use that 
\begin{equation}\label{a} v_p(n+m)=v_p(n) \quad \mbox{   if $v_p(n)<v_p(m)$;}
\end{equation} in particular, $v_p(n+p^x)=v_p(n)$ if $|n|<p^x$.

\begin{lemma}\label{4} 
Let $k$, $e$ and $p$ be as in Lemma \ref{lm:period}. Then
$$ v_p\left[(k+1){k\choose j}\right]\le e, \quad \mbox{for all $0\le j\le k$}.
$$
\end{lemma}
\begin{proof}
The proof hinges on the following fact: {\it Among all products of 
$x\le p^e(p-1)$ consecutive integers in the interval
$1,2,\dots,p^{e+1}-1$ the maximum valuation is attained by the
string that starts with $p^e$.} To show this, assume that
$r,r+1,\dots,r+x-1$ is such a string of integers. If $r>p^{e}$, say
$ap^{e}<r\le (a+1)p^{e}$, then the string beginning with $r-ap^{e}$
has the same valuation. Thus we may assume that $r\le p^{e}$.

If $r<p^{e}$, let $s$ be the least number such that $r<s\le p^{e}$
and $v_p(r)<  v_p(s)$. Then it is easy to see that
$v_p(r(r+1)\cdots (r+x-1))\le v_p(s(s+1)\cdots (s+x-1)$, and the
claim follows.

We may assume that $j\le \frac{k}{2}$. Then $j+1\le p^e(p-1)$,
and what was just shown implies that $v_p((k-j+1)(k-j+2)\cdots (k+1))\le
v_p(p^e(p^e+1)\cdots (p^e+j))= e+v_p(j!)$, which is equivalent to
the stated formula.
\end{proof}

\smallskip
\noindent {\sc Proof of Lemma \ref{lm:period}}  
{\ }We know from Pascal's triangle that ${d\choose k+1}\ 
\pmod {p^{r}}$ 
is completely determined by the values of ${j\choose 0}=1, j\ge 0,$ and 
${d'+i\choose i}\ \pmod{p^r}$ for any
$d'\ge 0$ and all $i=1,\dots,k+1$. Therefore it
suffices to show that 
\begin{equation}\label{eq1}
{p^{e+r}-1+i\choose i}\equiv_{p^r} 0,
\end{equation}
for all
$i=1,\dots,p^{e+1}-1,$ to establish the first part of the lemma. We have that
$v_p(p^{e+r}+s)=v_p(s)$ for all $s=1,\dots,p^{e+r}-1$. So the expansion of the
binomial coefficient
\[ {p^{e+r}-1+i\choose i}=\frac{(p^{e+r}+i-1)(p^{e+r}+i-2)\cdots p^{e+r}}
{i(i-1)(i-2)\cdots 2\cdot 1}
\] gives $v_p({p^{e+r}-1+i\choose i})=e+r-v_p(i)\ge r$, if $i<p^{e+1},$
which proves (\ref{eq1}).

The second part of the lemma is obvious when $0\le i\le k$, since both 
sides are zero (for the left-hand side this follows from $(i)$). Therefore,
assume that $k<i<p^{e+r}$. For each $j\neq 0$ write
$j=p^{\min \{v_p(j),e\}}q_j$. Then, for $0<j<p^{e+r}$:
\begin{equation}\label{q} q_{p^{e+r}-j}\equiv_{p^r}q_{-j} = -q_j.
\end{equation}

We have the equality 
\[
{i\choose k+1}\frac{(k+1)!}{p^{v_p((k+1)!)}} =
\frac{i(i-1)\cdots (i-k)}{p^{v_p((k+1)!)}} = \]
\[
p^{\sum_{j=0}^k \min\{v_p(i-j),e\}-v_p((k+1)!)}\prod_{j=0}^k q_{i-j} =
p^{\alpha}\prod_{j=0}^k q_{i-j},
\]
and similarly
\[ {p^{e+r}+k-i\choose k+1}\frac{(k+1)!}{p^{v_p((k+1)!)}} = 
\frac{(p^{e+r}-i)(p^{e+r}-i+1)\cdots
(p^{e+r}-i+k)}{p^{v_p((k+1)!)}} = \]
\[ p^{\sum_{j=0}^k
\min\{v_p(p^{e+r}-(i-j)),e\}-v_p((k+1)!)}\prod_{j=0}^k q_{p^{e+r}-(i-j)} =
p^{\beta}\prod_{j=0}^k q_{p^{e+r}-(i-j)}.
\]
We claim that $$\beta=\alpha\ge 0.$$
The equality follows from (\ref{a}), and the inequality will soon be 
proved. The two identities therefore give, using (\ref{q}):
$$
{i\choose k+1}\frac{(k+1)!}{p^{v_p((k+1)!)}} \equiv_{p^r} (-1)^{k+1}
{p^{e+r}+k-i\choose k+1}\frac{(k+1)!}{p^{v_p((k+1)!)}}. $$
Since  $\frac{(k+1)!}{p^{v_p((k+1)!)}}$ is invertible in $\Z_{p^r}$
this implies $(ii)$.

It remains to show that $\alpha\ge 0$. If $v_p(i-j)\le e$ for all
$j=0,\dots,k$, then $\alpha=v_p({i\choose k+1})\ge 0$. If not, then since
$k+1<p^{e+1}$ there is exactly one $s$, with $i-k\le s\le i$,
such that $v_p(s)>e$. In that case we have
$$
\alpha+v_p((k+1)!)=v_p((i-k)\cdots s\cdots i)-v_p(s)+e = $$
$$
v_p((i-k)\cdots (s-1))+v_p((s+1)\cdots i)+e =
v_p((s-(i-k))!)+v_p((i-s)!)+e,
$$ 
where the last equality uses (\ref{a}) twice.
Thus, using Lemma \ref{4},
$$
\alpha=-v_p((k-(i-s)+1)(k-(i-s)+2)\cdots (k+1))+v_p((i-s)!)+e = $$
$$
-v_p\left((k+1){k\choose i-s}\right)+e\ge 0.
$$

\smallskip
To prove $(iii)$  assume that ${d+i\choose k+1}\equiv_{p^r} 0$,
for some $d\ge k+1$ and all $i=0,\dots,k$.  
Then  ${d\choose j}\equiv_{p^r} 0$, for
$j=1,\dots,k+1$.  Especially, 
${d\choose p^s}\equiv_{p^r} 0$ for $0\le s\le e$,
which gives that $v_p({d\choose p^s})\ge r$. In particular,
$v_p(d)\ge r$.  We will now show that $v_p(d)\ge r+s$ for all $0\le s\le e$,
by induction on $s$.  Assume that $v_p(d)\ge r+s-1$. Then 
$r\le v_p({d\choose p^s})=v_p(d(d-1)\cdots(d-(p^s-1)))-
v_p(1\cdot 2\cdots(p^s-1)p^s)=v_p(d)+v_p((p^s-1)!)-v_p((p^s-1)!)-v_p(p^s)=
v_p(d)-s$. Hence, a run of $k+1$ consecutive zeros must begin with
${d\choose k+1}$ for some $d$ divisible by $p^{e+r}$. On the other
hand, the $1$'s along the left boundary of Pascal's triangle show that
there cannot be a run of more than $k+1$ zeros of the form ${i\choose k+1}$.
This proves the lemma.
\qed

\medskip
We can now proceed toward the proof of Theorem \ref{th:main}.
\begin{lemma}\label{lm:delare} For each 
$k<\left\lfloor\frac{d+1}{2}\right\rfloor$,
$G(d,k)$ is a divisor of 
$\frac{d-k+1}{\gcd(d-k+1, p_1^{e_1}p_2^{e_2}\cdots p_t^{e_t})}$, 
where $p_1,\dots,p_t$ are the primes $\le k+1$ and 
$p_i^{e_i}\le k+1<p_i^{e_i+1}$.
\end{lemma}

\begin{proof} Take a prime $p$ dividing $G(d,k)$ and let 
$x:=v_p(G(d,k))\ge 1$. 
Write
$k+1$ in base $p$, 
$k+1=k_0+k_1p+\cdots+k_ep^e$, where $0\le k_i<p$ and $k_e\neq 0$. Notice that
$$ v_p(\frac{d-k+1}{\gcd(d-k+1,p_1^{e_1}p_2^{e_2}\cdots p_t^{e_t})})
=\begin{cases} v_p(d-k+1)-e&\text{if $v_p(d-k+1)\ge e$,}\\
0&\text{otherwise,}
\end{cases}$$
so it will suffice to show that $v_p(d-k+1)-e \ge x$ in the first case and
obtain a contradiction in the second.

Since $p^x|G(d,k)$ we get that ${d+1-i\choose k+1}\equiv_{p^x}{i\choose
k+1}$, for all $i=1,\dots,\da$. Especially 
we must have ${d+1-i\choose k+1}\equiv_{p^x} 0$,
for $i=1,\dots,k$ and ${d-k\choose k+1}\equiv_{p^x} 1$.  From 
\[ (d-2k){d+1-k\choose k+1}= {d-k\choose k+1}(d-k+1),
\] we get $v_p(d-k+1)-v_p(d-2k)=v_p({d+1-k\choose k+1})\ge x\ge 1$.  Hence
by (\ref{a}),
$v_p(k+1)=v_p(d-k+1-(d-2k))=v_p(d-2k)<v_p(d-k+1)$.

There are now two cases:  First assume that $v_p(d-k+1)\ge e$. 
If $k+1=k_ep^e$ we are done, since we have shown that 
$v_p(d-k+1)-v_p(k+1)\ge x.$  Assume that $k+1>k_ep^e$.  From 
$d-k+1\ge k+1>k_ep^e$ we conclude that $v_p(d-k+1-k_ep^e)\ge e$, 
which implies $v_p(d-k+1-k_ep^e-i)=v_p(i)$, for all
$i=1,\dots,k+1-k_ep^e<p^e$. This in turn implies that 
\[v_p\left(\frac{(d-k-k_ep^e)!}{(d-2k)!}\right)=v_p((k-k_ep^e)!)=
v_p\left(\frac{(d+1-k_ep^e)!}{(d+1-k)!}\right).
\]
Using the equality
\[\frac{(d-k-k_ep^e)!}{(d-2k)!}{d+1-k_ep^e\choose k+1}={d+1-k\choose k+1}
\frac{(d+1-k_ep^e)!}{(d+1-k)!},
\] 
we get that $v_p({d+1-k_ep^e\choose k+1})=v_p({d+1-k\choose k+1})$.
This together with the identity
$$(d-k+1-k_ep^e){d+2-k_ep^e\choose k+1}=
{d+1-k_ep^e\choose k+1}(d+2-k_ep^e)$$  gives
\[x\le v_p(m_{k_ep^e-1,k})=v_p\left({d+2-k_ep^e\choose k+1}\right)\le \]
\[ v_p\left( {d+1-k_ep^e\choose
k+1}\right)-e+
v_p(d+2-k_ep^e)=
\]
\[v_p\left({d-k+1\choose k+1}\right)-e+v_p(k+1-k_ep^e)=v_p(d-k+1)-e.
\] 
Here the last equality comes from $v_p(k+1-k_ep^e)=v_p(k+1)=v_p(d-2k)$ 
and $v_p({d-k+1\choose k+1})=v_p(d-k+1)-v_p(d-2k)$, established above.

The second case is if $a:=v_p(d-k+1)<e$. The same argument can be applied 
again; however, now replacing $k_ep^e$ everywhere by
$k_ap^a+\cdots +k_ep^e$ and replacing $e$ by  $a$.  
We then get $x\le v_p(d-k+1)-a=0$, a contradiction.
\end{proof}

\begin{lemma}\label{lm:coroll}
$G(d,k)$ is a divisor to $m_{0,k}={d+1\choose k+1}$.
\end{lemma}
\newpage
\begin{proof} 
If for a prime $p$ we have that $p^r$ divides $G(d,k)$
and $p^e\le k+1 <p^{e+1}$, then  
Lemma \ref{lm:delare} implies that $p^{r+e}$
divides $d-k+1$. Hence $v_{p}({d+1\choose k+1})=v_{p}({d+1\choose k})
+v_{p}(d-k+1)-v_{p}(k+1)\ge r$. 
\end{proof}

\noindent {\sc Proof of Theorem \ref{th:main}.} The first statement follows 
from the fact that $m_{d-k,k}=1$ for 
$k\ge\left\lfloor\frac{d+1}{2}\right\rfloor.$

Let $k< \left\lfloor\frac{d+1}{2}\right\rfloor$. We have from the definition
of $m_{i,k}$ that for every prime $p$ and every $r\ge 1$:
\begin{equation}\label{c}
p^r|G(d,k)\Longleftrightarrow {d+1-i\choose k+1}\equiv_{p^r}{i\choose k+1},
\mbox{for $i=0,1,\dots,\left\lfloor\frac{d+1}{2}\right\rfloor.$}
\end{equation}
Actually, the definition supports this only for $i=1,\dots,
\delta= \left\lfloor\frac{d}{2}\right\rfloor$
on the right hand side, but $i=0$ can be added because of Lemma
\ref{lm:coroll} and $i= \left\lfloor\frac{d+1}{2}\right\rfloor $
(for $d$ odd) gives a trivially true identity.

{\em Case 1:} $k$ even. Assume that $G(d,k)\neq 1$, and that $p^r|G(d,k)$.
Since by Lemma \ref{lm:period}  there is a unique longest run of $k+1$
zeros in the  period of ${i\choose k+1}$ $\pmod {p^{e+r}}$ we get 
from (\ref{c}) that $d-k+1\equiv_{p^{e+r}} 0$. Therefore, Lemma
\ref{lm:period} and (\ref{c}) give
$$ {k+1\choose k+1}\equiv_{p^r}{d-k\choose k+1}
\equiv_{p^r}{p^{e+r}-1\choose k+1}\equiv_{p^r}-{k+1\choose k+1},
$$ which implies $p=2$ and $r=1$. Hence, $G(d,k)=2$, and this happens
only if $d-k+1\equiv_{2^{e+1}}0$. On the other hand, if
$d-k+1\equiv_{2^{e+1}}0$ then $2|G(d,k)$ can be concluded from
Lemma \ref{lm:period} and (\ref{c}).

{\em Case 2:} $k$ odd. Let $p^r$ be a divisor of 
$\frac{d-k+1}{\gcd(d-k+1, p_1^{e_1}p_2^{e_2}\cdots p_t^{e_t})}$. By
Lemma \ref{lm:delare} it suffices to show that $p^r|G(d,k)$.
The assumption implies that $p^{r+e}$ divides $d-k+1$, where as
usual $e$ is defined by $p^{e}\le k+1 <p^{e+1}$. Hence by
Lemma \ref{lm:period}
\[{d+1-i\choose k+1}\equiv_{p^r} {i\choose k+1},\quad
\mbox{ for all $i=0,\dots,d+1$},
\] which via (\ref{c}) shows that $p^r|G(d,k)$.

This  finishes the proof of the theorem.\qed

\medskip
\noindent
\begin{example} We want to calculate $G(116,9)$. Since $k=9$ is odd 
we calculate the greatest common divisor of $116-9+1=108$ and 
$2^3\cdot 3^2\cdot 5\cdot 7$ which is $36$.
We get $G(116,9)=118/36=3$.
 \qed
\end{example}

\section{Periodicity of $(d,k)$-realizable numbers}

We will now show the general theorem about the ultimately
stable periodic distribution of the $(d,k)$-realizable numbers.

\begin{theorem}\label{gap-thm} 
Fix $0\le k<d$, and let $G(d,k)$ be the number defined in equation
(\ref{Gdef}). Then there exists an integer $N$ such that for all 
$n>N$:

\[\mbox{$n$ is the number of $k$-faces of  
a simple $d$-polytope }\] 
\[\Longleftrightarrow\]
\[n\equiv 0\quad  \pmod {G(d,k)}\]
\end{theorem}
\begin{proof}
We will prove the theorem with the last statement replaced by 
$n\equiv m_{0,k}\; \pmod {G(d,k)}$. 
Lemma \ref{lm:coroll} shows that $m_{0,k}$ is 
divisible by $G(d,k)$, so this reformulation is equivalent.

\smallskip\noindent
$\Longrightarrow$   This direction is clear from Theorem \ref{g-thm}.

\smallskip\noindent
$\Longleftarrow$
Write
\[             G(d,k) = \sum_{i=1}^{\da} \lb_i m_{i,k},  \qquad \lb_i \in \Z.
\]
Suppose $m_{1,k} = C\cdot G(d,k)$. Define
\[       g_{\da}:=\left\{
\begin{array}{rl}
 (C-1) |\lb_{\da}| ,&\quad \mbox{ if $\lb_{\da} <0$}\\
 0 ,                &\quad\mbox{ otherwise;}
\end{array}\right.
\]
and recursively
\[       g_i := g_{i+1} + (C-1)(|\lb_{i}| + |\lb_{i+1}|) ,\qquad  0<i<\da.
\]
Let   $ N := m_{0,k}+ \sum_{i=1}^{\da} g_i m_{i,k},$
and let   $g_i^{(p)} = g_i + p \lb_i$, for    $p = 0,1,\dots$

Then  $g^{(p,q)} = (1, g_1^{(p)} + q, g_2^{(p)}, \dots , g_{\da}^{(p)})$
is nonnegative and decreasing after the first entry
for all $q \ge 0$ and all $0 \le p < C$, and hence is
an M-sequence. The $f_k$ values corresponding to these $g$-vectors are
\[
        f_k^{(p,q)} = N + qCG(d,k) + pG(d,k).     
\]
It is clear from the construction that all numbers $N+j\cdot G(d,k),
j=0,1,\dots$ are of the form $f_k^{(p,q)}$ for suitable
$q \ge 0$ and $0 \le p < C$.
\end{proof}

\begin{corollary} If the $m_{i,k}$ are relatively prime, then all numbers
from some point on are $(d,k)$-realizable. Furthermore, Theorem \ref{th:main}
shows that this happens precisely in the following cases:
\begin{enumerate} 
\item[(i)]{if $k\ge\left\lfloor\frac{d+1}{2}\right\rfloor$;}
\item[(ii)]{if $k<\left\lfloor\frac{d+1}{2}\right\rfloor$  is even, unless 
$d-k+1\equiv 0$ (mod $2^{e+1}$);}
\item[(iii)]{ if
$k<\left\lfloor\frac{d+1}{2}\right\rfloor$  is odd, unless
$d-k+1$ fails to divide $p_1^{e_1}p_2^{e_2}\cdots p_t^{e_t}.$}\qed
\end{enumerate}
\end{corollary}

Now, define $N(d,k)$ to be the least number $N$ for which Theorem 
\ref{gap-thm} is true. Note that 
$N(d,d-1)=d,$ so in what follows we may assume that $k<d-1$.

What can be said about the magnitude of $N(d,k)$? We will here give a general
upper bound, and then we will determine the exact
asymptotic growth for the special case $N(d,0)$ in the following section.

Define $$L(d,k):= \min \max_{i=1}^{\da} |\lb_i|,$$
with the minimum taken over all ways to represent $G(d,k)$ on the form
$$ G(d,k) = \sum_{i=1}^{\da} \lb_i m_{i,k},  \qquad \lb_i \in \Z.$$

\begin{lemma}
For all $0\le k\le d-2$, we have $L(d,k)< m_{1,k}$.
\end{lemma}
\begin{proof}
Assume $ G(d,k) = \sum_{i=1}^{\da} \lb_i m_{i,k},$ with 
$|\lb_s| \ge m_{1,k}$ and $m_{s,k}\neq 0$.  By symmetry we may 
assume that $\lb_s$ is positive, that is $\lb_s\ge m_{1,k}$.  
Since $G(d,k)<m_{1,k}m_{s,k}$, there has to be a $t$ such that $\lb_t<0$
and $m_{t,k}\neq 0$. Let 
\[\lb_i'=\begin{cases}\lb_s-m_{t,k}&\text{if $i=s$,}\\
\lb_t+m_{s,k}&\text{if $i=t$,}\\
\lb_i&\text{otherwise.}
\end{cases}
\]
We get $ G(d,k) = \sum_{i=1}^{\da} \lb_i' m_{i,k}.$ 
Since $|\lb_s'|<|\lb_s|$,
$|\lb_t'|<|\lb_t|$ or else $|\lb_t'|<m_{s,k}$,
and all other $|\lb_i'|$ are unchanged,  
we can continue this
process until $|\lb_i|<m_{1,k}$, for all $i$.
\end{proof}
\begin{theorem}\label{N-thm} $N(d,k) < \frac{1}{2}d^2 {d\choose k+1}^3$
\end{theorem}
\begin{proof}
Referring to the proof of Theorem \ref{gap-thm}, with an optimal
choice of the $\lb_i$-s, we have
$$N(d,k)\le \sum_{i=0}^{\da} g_i m_{i,k}\le {d+1\choose k+1}+
L(d,k) (C-1) \sum_{i=1}^{\da} (2\da +1-2i) {d+1-i\choose k+1}$$
$$\le {d+1\choose k+1}+ L(d,k) (C-1) \da (2\da -1) {d\choose k+1}<
2L(d,k) \da^2 {d\choose k+1}^2.$$
\end{proof}

For example, let $k=0$. The general
bound specializes to $N(d,0)\le \frac{1}{2}d^5.$
This should be compared to the true asymptotic value
$N(d,0)\sim cd^{3/2}$, which will be proved in the next section.

For $k\ge \left\lfloor\frac{d+1}{2}\right\rfloor$ 
we can improve on the general bound significantly.

\begin{theorem}\label{th:stora_k} 
Suppose $k\ge \left\lfloor\frac{d+1}{2}\right\rfloor$. Then
$N(d,k)< {d+1\choose d-k}(d-k)(k+1)(d+1)/2$. 
\end{theorem}

Since $G(d,k)=1$ for such $k$ this implies that
for every $n\ge {d+1\choose d-k}(d-k)(k+1)(d+1)/2$ there is a simple 
$d$-polytope with $n$ faces of dimension $k$.

To prove this we need a more technical construction than before.
First we extend the definition of $\pa^s$. 
Define for $p\le s$:
\[\pa_p^s(n):={a_s-p\choose s-p}+{a_{s-1}-p\choose s-1-p}+
\dots+{a_i-p\choose i-p},
\]
where $n$ is written in the unique expansion
\[ n={a_s\choose s}+{a_{s-1}\choose s-1}+\dots+{a_i\choose i},
\]
as in Section 2.  Also let $\pa_p^s(0):=0$. 
We will allow $p$ to be negative, which
corresponds to the natural ``inverse'' of $\pa_p^s$ for positive $p$.  Thus,
for $p>0$, $\pa^s_{-p}(n)$ is the greatest number such that 
$\pa^s_p(\pa^s_{-p}(n))=n$.
We will continue to write just $\pa^s$ for $\pa_1^s$.

Now, fix $d$ and $k\ge \left\lfloor\frac{d+1}{2}\right\rfloor$.  
Define a vector $\gv:=(g_0,g_1,\dots,g_{d-k})$
inductively as follows.  
\begin {itemize}
\item{} Let $g_{d-k}:=0$.
\item{} Assume we have defined $g_{d-k},g_{d-k-1},\dots,g_i$, for some 
$0<i\le d-k$.   Let 
\[g_{i-1}:=\pa^i(x_i),\quad \mbox{ where $x_i$ is the smallest integer
such that $x_i \ge g_i$  and}\]
\begin{equation}\label{g-vector}
\sum_{s=i}^{d-k}\left(\pa^i_{i-s}(x_i)-g_s \right)m_{s,k} 
\ge m_{i-1,k}-1.
\end{equation}
\end {itemize}
This is an $M$-sequence by construction.

\begin{lemma}\label{lm:constr}
Given the $g$-vector above, define
$N:=\sum_{i=0}^{d-k} g_i m_{i,k}$. Then there are no 
$(d,k)$-gaps larger than or equal to $N$.
\end{lemma}

\begin{proof}
Adding any positive integer to $g_1$ in an $M$-sequence gives another
$M$-sequence. Thus we only have to prove that it is possible to form all the
$m_{1,k}-1$ integers following $N$ with legal $g$-vectors.
This will imply the lemma.

We will think of the elements in column $k$ of $M_d$ as weights which
we combine to get the correct total weight.

Consider first the choice of $g_{d-k-1}:=\pa^{d-k}(x_{d-k})$, where 
$x_{d-k}= m_{d-k-1,k}-1$. All the vectors $(g_0,g_1,\dots,g_{d-k-1},i)$, 
$i=0,\dots m_{d-k-1,k}-1$ are $M$-sequences,
producing $N,N+1,\dots,N+m_{d-k-1,k}-1$ $k$-faces respectively.
We here use the fact that $m_{d-k,k}=1$. 
Similarly $(g_0,g_1,\dots,\allowbreak g_{d-k-1}+j,i)$, for fixed $j$ and
$i=0,\dots m_{d-k-1,k}-1$ gives $N+jm_{d-k-1,k},N+jm_{d-k-1,k}+1,
\dots,N+(j+1)m_{d-k-1,k}-1$ $k$-faces.  The definition of $g_{d-k-2}$ 
allows us to have $j$ sufficiently large to get all the numbers at least 
up to and including $N+m_{d-k-2,k}-1$.

Assuming inductively that we can form the sequence
$N, N+1, \dots,N+m_{i,k}-1$ by increasing only coordinates
$i+1,\dots,d-k$, the definition of $\gv$ gives that we 
can form all the numbers  $N, N+1,\dots,N+m_{i-1,k}-1$ 
by increasing only coordinates $i,\dots,d-k$ of $\gv$.
This proves the lemma
\end{proof}

\begin{example}
Take $d=10$ and $k=6$. We see from the matrix $M_{10}$, displayed
in Section $2$, that the
weights are $330,120,36,8$ and $1$.  We get $\gv=(1,4,6,6,0)$ which
gives $N(10,6)< 1074$, showing that every $n\ge 1074$ is 
$(10, 6)$-realizable.
\end{example}

\medskip
\noindent\begin{sc}Proof of Theorem \ref{th:stora_k}\end{sc}
First we show that $g_i\le(d-k-i)(k+1)$ by reverse induction.
It is trivially true for $g_{d-k}$.  Assume it is true for $g_i$.
Since $g_{i-1}=\pa^i(x_i)\le x_i$, it suffices to bound $x_i$.
The inequality \eqref{g-vector} is true if $(x_i-g_i)m_{i,k}\ge
m_{i-1,k}-1$.  Since $x_i$ is chosen to be minimal we get that
\[x_i\le g_i+\left\lceil\frac{m_{i-1,k}-1}{m_{i,k}} \right\rceil= 
g_i+\left\lceil\frac{{d+2-i\choose k+1}-1}{{d+1-i\choose k+1}} \right\rceil
\le\]
\[
\le g_i+\left\lceil\frac{d+2-i}{d+1-i-k}-\frac{1}{{d+1-i\choose k+1}}
\right\rceil\le g_i+k+1\le (k+1)(d-k-i+1).
\]
Now, 
\[N \le {d+1\choose k+1}+(k+1)\sum_{i=1}^{d-k}(d-k-i){d+1-i\choose k+1}\le\]
\[ {d+1\choose k+1}+(k+1){d\choose k+1}{d-k\choose 2}
\le {d+1\choose k+1}\left (1+\frac{(k+1)(d-k-1)(d+1)}{2}\right ).\]
This proves the theorem.
\qed

\medskip
Note that $m_{0,k}+1={d+1\choose d-k}+1$ is never $(d,k)$-
realizable for $k<d-1$.
This gives a trivial lower bound for $N(d,k)$ to be compared with the upper
bounds in Theorems  \ref{N-thm} and \ref{th:stora_k}.

\section{The case of vertices}\label{sc:vertices}
The only case of $(d,k)$-realizability that seems to have
been previously studied  is for $k=0$, i.e., the
number of vertices. We will make a more exact analysis of that case.

In \cite[Corollary 4.4.15]{Le}, Lee shows that for each dimension $d$
all sufficiently large numbers are $(d,0)$-realizable
(with parity restrictions, see the Introduction).  
Prabhu \cite {P1,P2} strengthened
the result and proved that there exists a constant $c$ such that every
$n>cd\sqrt{d}$ is $(d,0)$-realizable (with parity restrictions).  
This gives an upper bound on the
size $N(d,0)$ of the largest gap in each dimension --- we are not aware of any
published non-trivial lower bound.  The exact result is previously 
known only for small dimensions, see \cite {Le} where Lee 
lists all $(d,0)$-gaps for $d\le 9$.

We will sharpen Prabhu's result in both directions and prove that one can use
$c=\sqrt{2}+\epsilon$ as constant in his theorem for 
any $\epsilon>0$ and sufficiently large even $d$. 
However, the statement is not true for $c<\sqrt 2$.

\begin{theorem}\label{th:lower}
If $d\ge 4$ is even, then there does not exist a simple 
$d$-polytope with $n=(d-1)(\left\lceil \sqrt {2d-4}\right\rceil-2)+4$
vertices.\qed
\end{theorem}

\begin{theorem}\label{th:upper} For every even $d\ge 2$ and 
every $n>(d-1)(\sqrt{2d}+2\sqrt{2\sqrt{2d}}+5)$, there exists a 
simple $d$-polytope with $n$ vertices.\qed
\end{theorem}

\noindent
If we restrict our attention to $d$ odd, the true value for $c$ is
bounded by $1\le c\le 1+\epsilon$.

\begin{theorem}\label{th:odd}
If $d\ge 3$ is odd, then there does not exist a simple 
$d$-polytope with $n= (d-1)(\left\lceil \sqrt {d-2}\right\rceil-2)+4$
vertices, but for every even integer 
$n>(d-1)(\sqrt{d}+2\sqrt{2\sqrt{2d}}+5)$, there exists a 
simple $d$-polytope with $n$ vertices.\qed
\end{theorem}

\noindent
{\sc Proof of Theorem \ref{th:lower} }
Let $d=2\da\ge 6$, so the first column of $M_d$ will be $2\da+1,2\da-1,\dots,3,1$.
We will look for the lowest possible value for $f_0$ such that 
$f_0\equiv 4 \pmod {2\da-1}$.  The entries of the first column
will be the weights by which we seek to create the value of
$f_0$.  They are $2,0,-2,-4,\dots,-(2\da-2) \pmod {2\da-1}$. 
By the properties of $M$-sequences we have to take at least $k+2$
weights to obtain $f_0\equiv 4\equiv -2\da+5\pmod{2\da-1}$, where 
\[2+\sum_{i=0}^k -2i\le -(2\da-1)-(2\da-5),
\]
corresponding to the $M$-sequence $\gv=(1,1,\dots,1,1)$ with $k+2$ $1$'s.
This is equivalent to
\[k(k+1)\ge 4\da-4.
\]
Now, choose $k$ such that $k\le\sqrt{4\da-4}<k+1$. We then get that 
\[f_0\ge 2\da+1+\sum_{i=0}^k 2\da-1-2i=
2\da+1+(2\da-k-1)(k+1)>
\]
\[2\da+1+(2\da-1)\left\lceil \sqrt{4\da-4}\right\rceil-
(\left\lfloor \sqrt{4\da-4}\right\rfloor^2+\left\lfloor 
\sqrt{4\da-4}\right\rfloor)>
\]
\[
(d-1)(\left\lceil \sqrt{2d-4}\right\rceil-2)+4.
\]
Hence, $(d-1)(\left\lceil {\sqrt{2d-4}}\right\rceil-2)+4$ is a gap.

The result is easily seen to be true also for $d=4$.\qed

\medskip
\noindent
{\sc Proof of Theorem \ref{th:upper} }
Let $d=2\da\ge 4$ (the case $d=2$ is easily checked).
As above the first column of $M_d$ will be 
$2\da+1,2\da-1,\dots,3,1$.  First we note that if
$n+1,n+2,\dots,n+d-1$  are all realizable then every integer larger
then $n$ is realizable since we can just add $1$ to $g_1$ in the
corresponding $M$-sequences.

As in the previous proof we let $k_1$ be such that
$k_1\le\sqrt{4\da-4}<k_1+1$.
We consider the $M$-sequences $1=g_0=g_1=\dots=g_i$ and
$0=g_{i+1}=g_{i+2}=\dots$, for $0\le i \le k_1+1$.
The corresponding values for $f_0$ constitute one sequence of odd
residues and one sequence of even residues modulo $d-1$, with no distance
being larger than $2(k_1-1)$.
Now we choose $k_2$ such that  
\[\sum_{i=0}^{k_2} -2i\le -(2k_1-1)\Longleftrightarrow k_2+1>\sqrt{2(k_1-1)}.
\]
It is clear that the $M$-sequences $1=g_0,2=g_1=\dots=g_j$,
$1=g_{j+1}=g_{j+2}=\dots=g_i$ and
$0=g_{i+1}=g_{i+2}=\dots$, for $0\le j<i \le k_1+1$ and $j\le k_2$, give
values for $f_0 \pmod{d-1}$ where no residue is more than
$2(k_2-1)$ away from another residue of the same parity.
Continuing this process, we choose integers $k_1,k_2,\dots, k_s$,
as small as possible such that $\sqrt{2(k_{i-1}-1)}<k_i+1$, for $2\le
i\le s$.  We stop when we have reached $k_s=1$.  Hence, every possible
value for $f_0\!\pmod{d-1}$ can be obtained with an $M$-sequence that has
coordinates satisfying $g_i\le j$, whenever $k_{j+1}+1<i$.

So if $f_0$ is a gap then we must have
\[f_0<2\da+1+\sum_{i=0}^{k_1}(2\da-1-2i)+\sum_{i=0}^{k_2}(2\da-1-2i)+
\dots+\sum_{i=0}^{k_s}(2\da-1-2i)<
\]
\[2\da+1+(k_1+1+k_2+1+\dots+k_s+1)(2\da-1)<\mbox{[by induction]}
\]
\[2\da+1+(k_1+1+2(k_2+1))(2\da-1)<(d-1)(\sqrt{2d-4}+2\sqrt{2\sqrt{2d-4}-2}+5).
\]
This estimate suffices to show the theorem.  \qed

\medskip
\noindent
{\sc Proof of Theorem \ref{th:odd} }
The proof can be carried out in the same manner as the two previous
proofs.\qed

\end{document}